\documentclass{article}

\usepackage{graphicx}
\graphicspath{ {./images/} }% Required for inserting images
\usepackage{xcolor}
\usepackage{latexsym}
\usepackage{amsfonts}
\usepackage{enumerate}
\usepackage{multicol,caption}
\usepackage{graphicx}
\usepackage{amssymb}
\usepackage{amsmath}
\usepackage{epic}
\usepackage{amsthm}
\usepackage{graphicx,comment}
\usepackage{amscd}

\begin{document}

\parindent=0.2in
\parskip .1in

\color{blue}

\noindent {\Huge {\bf Engaging Students Through Math Competitions}}

\color{black}

\vspace*{.2in}

\noindent {\Large {\bf {\em B\'ela Bajnok}}\footnote{B\'ela Bajnok is a Professor of Mathematics at Gettysburg College and is the Director of the American Mathematics Competitions program of the MAA.  His email address is bbajnok@gettysburg.edu.}}

\vspace*{.4in}

\begin{multicols}{2}

\noindent
I was fortunate to grow up in Hungary, a country with a long and distinguished history in mathematical competitions.  Like many of my friends, I was keenly involved in these contests.  We liked training for them, being exposed to one beautiful problem after another; we enjoyed participating in them, even when we didn't end up winning; and we loved being in a community that was welcoming yet tight-knit.  I have remained involved with mathematical competitions ever since, and thus was excited and honored when, in 2017, I was asked to become the director of the American Mathematics Competitions (AMC) program of the Mathematical Association of America (MAA).  While the winners of our competitions -- especially those participating in the International Mathematical Olympiad (IMO) -- receive much recognition, our goals are to encourage students to engage in mathematics nationwide, and to discover, develop, and nurture talent more broadly.

With the participation of over 300,000 students each year, the AMC is the largest program of the MAA.  The competitions start with the  AMC 8, the  AMC 10, and the  AMC 12 exams, open to students in grade 8 or below, grade 10 or below, and grade 12 or below, respectively.  Based on their performance on these multiple-choice competitions, approximately 6,000 students are invited to take the  American Invitational Mathematics Exam (AIME), a challenging three-hour exam where the answer to each of the 15 problems is a nonnegative integer under 1000.  The competition series culminates with the  USA Mathematical Olympiad (USAMO) and the  USA Junior Mathematical Olympiad (USAJMO), offered to approximately 500 students.  These Olympiads follow the style of the  IMO: they ask students to provide rigorous proofs for three problems on each of two consecutive days, with an allowed time of four and a half hours each day.  We hope that our competitions show students -- and their teachers and society at large -- the beauty and power of mathematics, usually well beyond what they see in typical mathematics classrooms.   

The creation of problems suitable for the AMC is a highly challenging task.  
Even at the beginning levels where problems stay close to the standard school curriculum, we aim for problems that ask questions in novel and fun ways.  Coming up with beautiful -- yet still elementary -- problems at the Olympiad level is particularly challenging.  Yet, year after year, collections of ingenious and captivating problems are created; everyone with an interest is encouraged to explore them.  Beyond the competitions themselves, we hope that our problems provide learning opportunities for current and future students, their teachers, and anyone with an interest in mathematics.  We present three examples here: one from an AMC 8 exam, one from an AIME competition, and one from a USA Mathematical Olympiad.

{\bf 2022 AMC 8 Problem 25} (created by Silva Chang):  
A cricket hops between four leaves, on each turn hopping to one of the other three leaves with equal probability. After four hops what is the probability that the cricket returns to the leaf where it started?

    \begin{center}
        \includegraphics[scale=0.7]{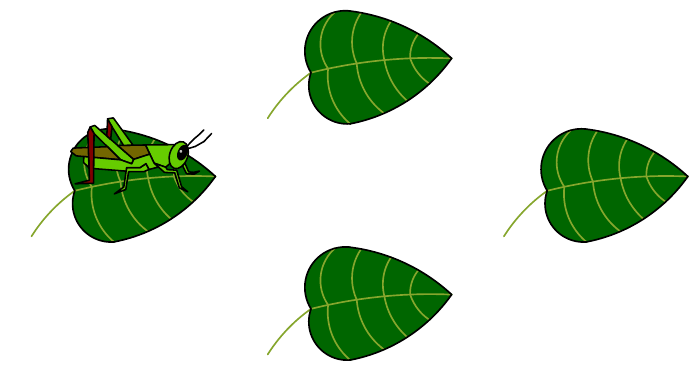}
        \captionof{figure}{Four leaves and a cricket}
    \end{center}

There is a large variety of different techniques to solve this problem: the official solution guide lists four of them, from case counting to recursions, and additional methods such as dynamic programming or generating functions can be used as well.  Problems being suitable for multiple approaches are of course helpful for the contestants at the time, but this subsequently also allows teachers to discuss techniques that are new to their students.

{\bf 2022 AIME I Problem 1} (created by David Altizio): Quadratic polynomials $P(x)$ and $Q(x)$ have leading coefficients $2$ and $-2,$ respectively. The graphs of both polynomials pass through the two points $(16,54)$ and $(20,53).$ Find $P(0) + Q(0).$

This problem can, of course, be solved by finding the two quadratic polynomials that satisfy the conditions of this problem, but this would take a while, and the answers are not particularly nice.  The clever approach (perhaps suggested by the question itself) would focus on the polynomial $P(x)+Q(x)$, which is linear; since its graph goes through the points $(16,108)$ and $(20,106)$, it equals $-x/2+116$, and thus $P(0)+Q(0)=116$.

One of the challenges we have when assembling our exams is to be able to select the students who will then move on to subsequent competitions and win awards and prizes, while also providing approachable problems to all participants.  Three hours may seem long for the AIME competition, but time is an important factor for many participants, so being able to solve the easiest problems quickly is a big advantage.     

{\bf 2021 USAMO Problem 3} (created by Shaunak Kishore and Alex Zhai):
Let $n \ge 2$ be an integer.
	An $n$-by-$n$ board is initially empty.
	Each minute, you may perform one of two moves:
\newline
	\hspace*{.1in} \textbullet\ If there is an L-shaped tromino region
		of three cells without stones on the board
		(see Figure 2; rotations not allowed),
		you may place a stone in each of those cells.
\begin{center}
	\includegraphics{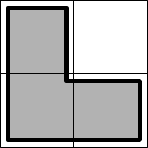}
 \captionof{figure}{The tromino region}
\end{center}	
	\hspace*{.1in} \textbullet\ If all cells in a row or column have a stone, you may remove all stones from that row or column.
 \newline
For which $n$ is it possible that, after some nonzero number of
	moves, the board has no stones?

This problem proved to be one of the most challenging ones in the history of the USAMO: only seven students who took the exam were able to solve it.  The solution employs what has recently been called the polynomial method, which is based on the simple fact that polynomials have a limited number of roots.  (Note that one can prove this fact with some basic algebra.) This technique has seen a lot of interesting applications lately in both mathematics and computer science, and we were excited to feature it in this beautiful problem.

It is not hard to see that when $n$ is divisible by $3$, the empty board is achievable. Divide the board into 3-by-3 sub-boards.  In each 3-by-3 sub-board, follow the procedure shown in Figure 3.

\begin{center}
\setlength{\unitlength}{.7cm}
\begin{picture}(11,3)

\put(0,0){\line(1,0){3}}
\put(0,1){\line(1,0){3}}
\put(0,2){\line(1,0){3}}
\put(0,3){\line(1,0){3}}

\put(0,0){\line(0,1){3}}
\put(1,0){\line(0,1){3}}
\put(2,0){\line(0,1){3}}
\put(3,0){\line(0,1){3}}

\put(4,0){\line(1,0){3}}
\put(4,1){\line(1,0){3}}
\put(4,2){\line(1,0){3}}
\put(4,3){\line(1,0){3}}

\put(4,0){\line(0,1){3}}
\put(5,0){\line(0,1){3}}
\put(6,0){\line(0,1){3}}
\put(7,0){\line(0,1){3}}

\put(4.5,.5){\circle*{0.5}}
\put(5.5,.5){\circle*{0.5}}
\put(4.5,1.5){\circle*{0.5}}

\put(8,0){\line(1,0){3}}
\put(8,1){\line(1,0){3}}
\put(8,2){\line(1,0){3}}
\put(8,3){\line(1,0){3}}

\put(8,0){\line(0,1){3}}
\put(9,0){\line(0,1){3}}
\put(10,0){\line(0,1){3}}
\put(11,0){\line(0,1){3}}

\put(8.5,.5){\circle*{0.5}}
\put(9.5,.5){\circle*{0.5}}
\put(8.5,1.5){\circle*{0.5}}
\put(9.5,2.5){\color{blue} \circle*{0.5}}
\put(9.5,1.5){\color{blue} \circle*{0.5}}
\put(10.5,1.5){\color{blue} \circle*{0.5}}

\end{picture}
\end{center}

\begin{center}
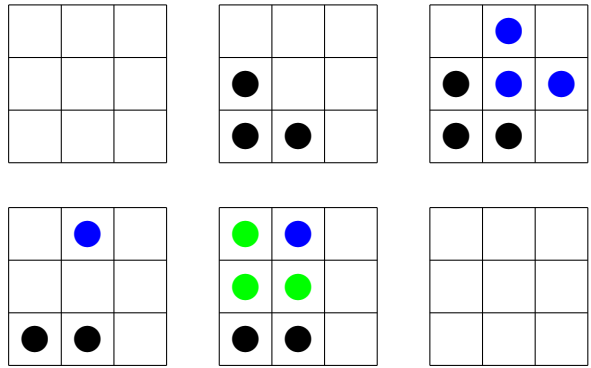

\setlength{\unitlength}{.7cm}
\begin{picture}(11,3)

\put(0,0){\line(1,0){3}}
\put(0,1){\line(1,0){3}}
\put(0,2){\line(1,0){3}}
\put(0,3){\line(1,0){3}}

\put(0,0){\line(0,1){3}}
\put(1,0){\line(0,1){3}}
\put(2,0){\line(0,1){3}}
\put(3,0){\line(0,1){3}}

\put(0.5,.5){\circle*{0.5}}
\put(1.5,0.5){\circle*{0.5}}
\put(1.5,2.5){\color{blue} \circle*{0.5}}

\put(4,0){\line(1,0){3}}
\put(4,1){\line(1,0){3}}
\put(4,2){\line(1,0){3}}
\put(4,3){\line(1,0){3}}

\put(4,0){\line(0,1){3}}
\put(5,0){\line(0,1){3}}
\put(6,0){\line(0,1){3}}
\put(7,0){\line(0,1){3}}

\put(4.5,.5){\circle*{0.5}}
\put(5.5,.5){\circle*{0.5}}
\put(5.5,2.5){\color{blue} \circle*{0.5}}
\put(4.5,2.5){\color{green} \circle*{0.5}}
\put(5.5,1.5){\color{green} \circle*{0.5}}
\put(4.5,1.5){\color{green} \circle*{0.5}}

\put(8,0){\line(1,0){3}}
\put(8,1){\line(1,0){3}}
\put(8,2){\line(1,0){3}}
\put(8,3){\line(1,0){3}}

\put(8,0){\line(0,1){3}}
\put(9,0){\line(0,1){3}}
\put(10,0){\line(0,1){3}}
\put(11,0){\line(0,1){3}}

\end{picture}
\captionof{figure}{Clearing the board when $n$ is divisible by 3}
\end{center}

%%%%%%%%%%%%%%%%%%%%%%%%%%%%%%%%%%%%%%%%%%%%%%%%%%%%%%%%%%%%%%%%%%%%%%%%%%%%

We can then use the polynomial method to prove that no procedure will succeed when $n$ is not divisible by $3$.  Leaving a bit of mystery, we only show this here for $n=4$; readers are encouraged to think about how this generalizes to other $n$ that is not divisible by $3$ (and why the argument fails when $n$ is divisible by $3$).   

We place monomials into the $16$ squares of our 4-by-4 board as shown in Figure 4.

\begin{center}
\setlength{\unitlength}{1cm}
\begin{picture}(5,5)
\put(0,1){\line(1,0){4}}
\put(0,2){\line(1,0){4}}
\put(0,3){\line(1,0){4}}
\put(0,4){\line(1,0){4}}
\put(0,5){\line(1,0){4}}

\put(0,1){\line(0,1){4}}
\put(1,1){\line(0,1){4}}
\put(2,1){\line(0,1){4}}
\put(3,1){\line(0,1){4}}
\put(4,1){\line(0,1){4}}

\put(.4,.5){0}
\put(1.4,.5){1}
\put(2.4,.5){2}
\put(3.4,.5){3}

\put(-.6,1.4){0}
\put(-.6,2.4){1}
\put(-.6,3.4){2}
\put(-.6,4.4){3}

\put(.4,1.4){1}
\put(.4,2.4){$y$}
\put(.4,3.4){$y^2$}
\put(.4,4.4){$y^3$}

\put(1.4,1.4){$x$}
\put(1.3,2.4){$xy$}
\put(1.3,3.4){$xy^2$}
\put(1.3,4.4){$xy^3$}

\put(2.3,1.4){$x^2$}
\put(2.2,2.4){$x^2y$}
\put(2.1,3.4){$x^2y^2$}
\put(2.1,4.4){$x^2y^3$}

\put(3.3,1.4){$x^3$}
\put(3.2,2.4){$x^3y$}
\put(3.1,3.4){$x^3y^2$}
\put(3.1,4.4){$x^3y^3$}

\end{picture}
\captionof{figure}{Labelling cells of a 4-by-4 board}
\end{center}

We keep track of our moves as follows: when stones are placed on the board, we add their values, and when stones are removed from the board, we subtract their values.  Assume now, indirectly, that a procedure results in an empty board; to be specific, suppose that in this procedure, a tromino with its lower left corner in position $(i,j)$ was added $t_{i,j}$ times; the $i$th column was cleared $c_i$ times; and the $j$th row was cleared $r_j$ times.  This then means that we must have
$$T(x,y)-C(x,y)-R(x,y)=0,$$
where
$$T(x,y)=\sum_{i=0}^2 \sum_{j=0}^2 t_{i,j}x^iy^j(1+x+y),$$
$$C(x,y)=\sum_{i=0}^3 x^i (1+y+y^2+y^3),$$
$$R(x,y)=\sum_{j=0}^3 y^j (1+x+x^2+x^3).$$

Consider now the set $D=\{-1, \mathrm{i}, - \mathrm{i}\}$, and observe that $C(x,y)$ and $R(x,y)$ both equal zero for every $x \in D$ and $y \in D$, but $1+x+y$ is never equal to zero for any such $x$ and $y$.  We then must have 
$$P(x,y)=\sum_{i=0}^2 \sum_{j=0}^2 t_{i,j}x^iy^j=0$$ for each $x \in D$ and $y \in D$.  This is a contradiction, however, since the polynomial $P(x,y)$ is at most quadratic in both $x$ and $y$, and thus it cannot have more than two values of each variable for which it equals zero.

%%%%%%%%%%%%%%%%%%%%%
   
The AMC program would not be possible without our four remarkable editorial boards that create our exams each year: they propose problems, review and rate all submissions, select the problems for the exams, and meticulously edit these problems and their solutions.  They also are instrumental in helping us shape general AMC policies and practices. These nearly 180 individuals come from all across the nation and from abroad; include people from academia, high schools, business, and industry; have a wide range of ages represented, from undergraduates to retired mathematicians; and possess an impressive variety of mathematical and cultural expertise.  I am  immensely grateful to this talented, diverse, and accomplished group of people.  Anyone with an interest in joining us is encouraged to contact me.  

{\small Credits}

Figure 1 is courtesy of Ricardo Concei\c{c}\~ao.

Figure 2 is courtesy of Shaunak Kishore and Alex Zhai.

Figures 3 and 4 are courtesy of the author.

Photo of B\'ela Bajnok is courtsey of Steve Berg.

\vspace*{.2in}

\begin{flushright}
\noindent {\small DRAFT September 2, 2023}
\end{flushright}

\end{multicols}

\end{document}